# AN EXAMINATION OF THE NEGATIVE OCCUPANCY DISTRIBUTION AND THE COUPON-COLLECTOR DISTRIBUTION

BEN O'NEILL[*], *Australian National University*[**]

WRITTEN 1 OCTOBER 2021

**Abstract**

We examine the negative occupancy distribution and the coupon-collector distribution, both of which arise as distributions relating to hitting times in the extended occupancy problem. These distributions constitute a full solution to a generalised version of the coupon collector problem, by describing the behaviour of the number of items we need to collect to obtain a full collection or a partial collection of any size. We examine the properties of these distributions and show how they can be computed and approximated. We give some practical guidance on the feasibility of computing large blocks of values from the distributions, and when approximation is required.

EXTENDED OCCUPANCY PROBLEM; NEGATIVE OCCUPANCY DISTRIBUTION; COUPON COLLECTOR DISTRIBUTION; COUPON COLLECTOR PROBLEM.

The "coupon-collector problem" (also sometimes known as the "dixie cup problem") involves a situation where we have an exchangeable sequence of draws with replacement from a finite set of objects, and our goal is to collect a full set of these objects. The problem can be traced back at least to Feller (1950, p. 213) and the asymptotics of the problem for large sets of objects was considered in Brayton (1963) and Baum and Billingsley (1965). Outlines of the problem and its history are found in Johnson and Kotz (1977, Section 3.5) and Holst (1986), and the problem has been revisited in a number of later papers (Dawkins 1991; Levin 1992; Read 1998; Kuenon 2000; Adler, Oren and Ross 2003; and Doumas and Papanicolaou 2012).

All versions of the coupon collector problem arise in the context of sampling with replacement. In the classical version of this problem the sampling distribution is uniform (i.e., each object is equally likely to be selected on each draw), though this can be generalised to allow for non-uniform sampling. In either case, the number of objects that have to be sampled to get a full collection is a random variable, and the goal of the problem is to describe some aspect of the behaviour of this random variable. Stipulations of the problem commonly ask for the expected number of objects sampled to obtain a full collection, but it is possible to obtain the distribution, which gives a complete specification of the behaviour.

---

[*] E-mail address: ben.oneill@anu.edu.au; ben.oneill@hotmail.com.
[**] Research School of Population Health, Australian National University, Canberra ACT 0200.

The coupon-collector problem can be generalised in several ways, and in this paper we consider two important generalisations that arise in the "extended occupancy problem". In that problem, a finite number of balls are randomly allocated bins (via a uniform distribution) in a manner that is analogous to sampling with replacement. In addition, there is an allowance for each ball to "fall through" its allocated bin with some fixed probability (see e.g., Samuel-Cahn 1974). We liken this to the coupon-collector process, where a ball falling through its bin corresponds to a "damaged" object that does not count towards the collection. The problem is simplest to analyse by looking at the stochastic process for the "occupancy number", which is the number of occupied bins after any given number of balls have been allocated. O'Neill (2021a) examines the occupancy number as a "pure birth" Markov chain, and uses this to derive the distribution of the hitting times for each fixed value of the occupancy number. This hitting time is the number of objects required to obtain a collection of any given size, which includes the number of objects required for a full collection.

We now frame the problem formally, using a stochastic process for the extended occupancy problem. Suppose we sequentially allocate balls (uniformly) at random to $m$ bins. Each ball has fixed probability $\theta$ of occupying its allocated bin and corresponding probability $1 - \theta$ of "falling through" that bin. Let $K_n$ denote the number of bins occupied by the first $n$ balls (called the "occupancy number") and note that $0 \leq K_n \leq \min(n, m)$. The stochastic process $\{K_n | n \in \mathbb{N}\}$ is "pure birth" process, which is a type of Markov chain; each allocated ball either increases the occupancy number by one, or it makes no change to the occupancy number. To obtain an occupancy number of $k$ bins (i.e., having $k$ bins occupied) we must allocate at least $k$ balls. It is useful to examine the *excess* number of balls required to occupy $k$ bins, which we denote as $T_k \equiv \min\{n \in \mathbb{N} | K_{k+n} = k\}$. O'Neill (2021a) shows that this random variable has a negative occupancy distribution $T_k \sim \text{NegOcc}(m, k, \theta)$ described by the mass function below. In the coupon-collector problem we look for the number of objects required for a full collection with $K_n = m$, and in this special case we have $T_m \sim \text{CoupColl}(m, \theta)$ described below.

**DEFINITION (The negative occupancy distribution):** This is a discrete distribution that has probability mass function given over all integer arguments $t \geq 0$ as:

$$\text{NegOcc}(t|m, k, \theta) = \frac{\theta^{k+t}}{m^{k+t}} \cdot (m)_k \cdot S\left(k + t - 1, k - 1, m \cdot \frac{1 - \theta}{\theta}\right),$$

where $0 < k \leq m \leq \infty$ are the **occupancy parameter** (number of occupied bins) and **space parameter** (number of bins) respectively, and $0 < \theta \leq 1$ is the **probability parameter**. □

**DEFINITION (Coupon-collector distribution):** This is a discrete distribution with probability mass function given over all integer argument values $t \geq 0$ as:

$$\text{CoupColl}(t|m, \theta) = \frac{\theta^{m+t}}{m^{m+t}} \cdot m! \cdot S\left(m + t - 1, m - 1, m \cdot \frac{1-\theta}{\theta}\right).$$

where $0 < m < \infty$ is the **space parameter** (number of bins) and $0 < \theta \leq 1$ is the **probability parameter**. Since $\text{CoupColl}(t|m, \theta) = \text{NegOcc}(t|m, m, \theta)$ the coupon-collector distribution is a special case of the negative occupancy distribution with $k = m$. □

In the above expressions, the terms $(m)_k = \prod_{i=0}^{k-1}(m-i)$ are falling factorials and the terms $S(n, k, \phi)$ are the noncentral Stirling numbers of the second kind, where the last term is the noncentrality parameter (for details of various properties of these numbers, see Charalambides 2005, pp. 73-96; O'Neill 2021a, Appendix I). The noncentral Stirling numbers are generally computed using recursive formulae, and so the mass function for the negative occupancy distribution is computationally onerous for large values of $k$ and $t$.

The coupon-collector distribution solves a generalised version of the classical coupon-collector problem by giving a behavioural description of the stochastic behaviour of the excess number of objects we need to collect to obtain a full collection. (Note that the distribution is for the excess number of objects collected; to obtain the distribution of the total number of objects collected we simply shift the location of this distribution $m$ units.) The parameter $\theta$ generalises the classical problem by giving a fixed probability that a given object is "damaged" and so does not contribute to the collection. The negative occupancy distribution further generalises things by allowing us to describe the stochastic behaviour of the excess number of objects we need to collect to obtain a partial collection of any given size.

In this paper, we extend O'Neill (2021a) by examining the properties of the negative occupancy distribution and the coupon-collector distribution in greater detail. The latter distribution is a special case of the former, but it is important to note up front that by "folding in" the parameter $k = m$ it has some different properties with respect to the (now dual) parameter $m$, particular with regard to its recursive properties. We will look at the probability functions and moments of the negative occupancy distribution, and various recursive equations and representations that relate the distribution to well-known distributional families. We also examine computation of the mass function of the negative occupancy distribution and see how it can be approximated in cases where exact computation is infeasible.

# 1. The negative occupancy distribution is a convolution of geometric distributions

The negative occupancy distribution can be derived by various methods from the underlying Markov chain describing the extended occupancy process. O'Neill (2021a) gives a derivation directly in terms of the underlying distribution of the occupancy number. However, another useful representation can be obtained by considering the total excess balls as a sum of the increments measuring the excess balls required for each unit increase in the occupancy number. For all $k = 1, 2, \ldots, m$ we define the increments $S_k \equiv T_k - T_{k-1}$ (with $T_0 \equiv 0$) as the excess number of balls required to move to the occupancy number $k$ from the previous value. Since the occupancy process is a Markov chain, the increments $S_1, S_2, \ldots, S_m$ are independent, and since these random variables each count the number of "failures" before a single "success", each follows a geometric distribution:[1]

$$S_\ell \sim \text{Geom}\left(\theta \cdot \frac{m - \ell + 1}{m}\right).$$

This distribution reflects the fact that progression of the occupancy number by a single unit requires a single ball to occupy a bin that is not already occupied. The probability of "success" is equal to the probability that the ball occupies its allocated bin (i.e., does not "fall through" that bin), multiplied by the probability that the ball occupies one of the $m - k + 1$ bins that was not already occupied.

Since $T_k = \sum_{\ell=1}^{k} S_\ell$ we see that a negative occupancy random variable is a sum of independent (but not identical) geometric random variables, and so the negative occupancy distribution must be a convolution of geometric distributions. We formalise this in Theorem 1 below (proofs of all theorems are in the Appendix).

**THEOREM 1 (Convolution representation):** The negative occupancy distribution satisfies the recursive convolution equation:

$$\text{NegOcc}(m, k, \theta) = \text{NegOcc}(m, k - 1, \theta) * \text{Geom}\left(\theta \cdot \frac{m - k + 1}{m}\right),$$

with the base case $\text{NegOcc}(m, 1, \theta) = \text{Geom}(\theta)$. This gives the convolution representation:

$$\text{NegOcc}(m, k, \theta) = \text{Geom}(\theta) * \text{Geom}\left(\theta \cdot \frac{m - 1}{m}\right) * \ldots * \text{Geom}\left(\theta \cdot \frac{m - k + 1}{m}\right).$$

---
[1] Here we refer to the geometric distribution on the support $s = 0, 1, 2, \ldots$ (i.e., counting only the "failures" not the total trials), so the mass function is $\text{Geom}(s|p) = (1 - p)^s p$.

Taking $m = \infty$ gives $S_\ell \sim \text{Geom}(\theta)$ for all $\ell = 1,2,3,...$ which means the above convolution reduces down to the $k$-fold convolution of a geometric distribution with parameter $\theta$, which is well-known to be the negative binomial distribution. Thus, we see that the negative occupancy distribution is a generalisation of the negative binomial distribution, which occurs in the special case where $m = \infty$. (This same limit result is established directly from the negative occupancy mass function in O'Neill 2021a, Theorem 7). Observe that we allow an infinite number of bins in the parameterisation of the negative occupancy distribution, but not for the coupon-collector distribution. In the latter case, we require $m < \infty$ to ensure that the collection will eventually terminate, otherwise we never obtain a full collection.[2]

The negative occupancy distribution is a broad family of discrete distributions that generalises the negative binomial distribution in a way that allows the user to describe the behaviour of excess hitting times in the extended occupancy problem. The imposition of a finite number of bins "squashes" the occupancy number and so the resulting hitting times tend to be larger than in the negative binomial case. In addition to being a generalisation of the negative binomial, the distribution also constitutes a special case of the broader class of distributions that are convolutions of geometric distributions. Sen and Balakrishnan (1999) give a general form for the latter class of distributions in cases where the probabilities parameters for the geometric distributions are distinct. This general form gives the probability mass function as a weighted sum of geometric distributions, leading us to the following representation.

**THEOREM 2 (Weighted geometric representation):** The negative occupancy distribution can be written as a weighted sum of geometric distributions:

$$\text{NegOcc}(t|m,k,\theta) = \sum_{i=1}^{k} w_{i,k} \cdot \text{Geom}\left(t + k - 1 \middle| \theta \cdot \frac{m-i+1}{m}\right),$$

using the weights:

$$w_{i,k} = (-1)^{k-i} \frac{(m)_k}{(m-i+1)(i-1)!(k-i)!}.$$

---

[2] Another way to put this is to observe that, for all $0 < \theta < 1$, as $m \to \infty$ we have:

$$\frac{\text{CoupColl}(t+1|m,\theta)}{\text{CoupColl}(t|m,\theta)} \to (1-\theta) \cdot \frac{m+t}{t+1} \to \infty.$$

(This result uses a limiting form in Lemma 4 of O'Neill 2021.) This means that as $m \to \infty$ the coupon-collector distribution moves (in an informal sense) towards a point-mass distribution on $t = \infty$. This reflects the fact that with an infinite number of bins, the excess time required to occupy all the bins must be infinite.

**REMARK:** The weights in this representation can be computed rapidly using recursive methods, using the following recursive form. For indices $1 \leq i < k \leq m$ we have:

$$\frac{w_{i+1,k}}{w_{i,k}} = -\frac{m-i+1}{m-i} \cdot \frac{k-i}{i} \qquad w_{k,k} = \frac{(m)_{k-1}}{(k-1)!}.$$

We can see that the leading weight $w_{k,k}$ is positive and the remaining weights alternate in sign. Since some weights are negative, this means that the representation is not an implicit mixture distribution (i.e., the weights cannot be considered as a probability mass function over the index $i$). We also see that the weights depend on the parameters $k$ and $m$ but not the probability parameter $\theta$, so they can be computed without consideration for this latter parameter. $\square$

The above exposition looks at the marginal case where we are concerned with the excess hitting time from the start of the occupancy process (where the occupancy number is zero). However, the family of negative occupancy distributions is "closed under conditioning" in the sense that it also gives the distribution of the excess hitting times between *any* two occupancy numbers. To see this, let $T_{k|r} \equiv T_{r+k} - T_r$ denote the excess hitting time between the occupancy numbers $r$ and $r+k$ (with $0 \leq r < r+k \leq m$). Since $T_{k|r} = \sum_{\ell=r+1}^{r+k} S_\ell$ we can write the distribution of this quantity as a convolution of geometric distributions, and this simplifies to:

$$T_{k|r} \sim \text{NegOcc}\left(m-r, k, \theta \cdot \frac{m-r}{m}\right).$$

This conditional variation of the problem includes the marginal problem as the special case where $r = 0$. It is therefore useful to give the generalised probability mass values:

$$\begin{aligned}
\mathbb{P}(T_{k|r} = t) &= \mathbb{P}(T_{r+k} - T_r = t) \\
&= \text{NegOcc}\left(t \,\middle|\, m-r, k, \theta \cdot \frac{m-r}{m}\right) \\
&= \frac{\theta^{k+t}}{m^{k+t}} \cdot (m-r)_k \cdot S\left(k+t-1, k-1, r+m \cdot \frac{1-\theta}{\theta}\right).
\end{aligned}$$

This equation gives a general solution to the behaviour of the excess hitting times between any two values of the occupancy number. In the case where $k = m$ this reduces to the coupon-collector distribution (with the same size and probability parameters), which can be used in conditional versions of the coupon-collector problem. The conditional distribution can easily be established formally using the convolution representation. However, it also has a simple intuition that elucidates the result. Since $r$ bins are already occupied, the probability of a ball being allocated to one of those pre-occupied bins can be "folded into" the probability parameter (giving the adjusted probability parameter above) and the number of bins can then be reduced to exclude those occupied bins (reducing the space parameter from $m$ to $m-r$).

The conditional result shown above demonstrates why the extended occupancy problem is such a useful extension of the classical occupancy problem. By incorporating a probability that the ball will "fall through" its bin, we can obtain distributions for conditional problems by taking the conditional occupancy number and "folding it" into the probability parameter. This allows us to solve versions of the coupon-collector problem (and broader hitting-time problems) where some coupons are already collected. In such cases, we can condition on the number of coupons already collected, and the coupon-collector distribution (or negative occupancy distribution) is then used to describe the behaviour of the number of remaining coupons to collect to collect a full (or partial) set.

## 2. Generating functions and moments of the negative occupancy distribution

Our previous section shows some basic representations that elucidate the negative occupancy distribution. In this section we derive the generating functions and moments of the distribution. The distribution has a closed-form expression for its cumulant function, and the derivatives of this function, which means we obtain closed-form expressions for all of its moments. We show the generating functions and the resulting cumulants for the distribution below.

**THEOREM 3 (Generating functions):** The negative occupancy distribution has the following generating functions. The probability generating function $G$, characteristic function $\phi$, moment generating function $M$, and cumulant function $K$ are given respectively by:

$$G(z) = \theta^k \cdot \prod_{\ell=m-k+1}^{m} \frac{\ell}{m - (m - \ell\theta)z} \qquad |z| < \frac{m}{m - (m-k+1)\theta},$$

$$\phi(s) = \theta^k \cdot \prod_{\ell=m-k+1}^{m} \frac{\ell}{m - (m - \ell\theta)e^{is}} \qquad |s| < \ln\left(\frac{m}{m - (m-k+1)\theta}\right),$$

$$M(s) = \theta^k \cdot \prod_{\ell=m-k+1}^{m} \frac{\ell}{m - (m - \ell\theta)e^{s}} \qquad s < \ln\left(\frac{m}{m - (m-k+1)\theta}\right),$$

$$K(s) = k \ln(\theta) + \sum_{\ell=m-k+1}^{m} \ln\left(\frac{\ell}{m - (m - \ell\theta)e^{s}}\right) \qquad s < \ln\left(\frac{m}{m - (m-k+1)\theta}\right).$$

(Note that the inequality for the domain of convergence for each expression ensures that each denominator term is greater than zero, so that the expression given is well-defined.)

The cumulants and moments of the negative occupancy distribution involve the generalised harmonic numbers. To facilitate this analysis, we denote these numbers by $H_m^{(r)} \equiv \sum_{\ell=1}^{m} \ell^{-r}$ (see e.g., Spieß 1990) and we define the function:

$$h_i(m, k, \theta) \equiv \left(\frac{m}{\theta}\right)^i (H_m^{(i)} - H_{m-k}^{(i)}).$$

For all $i = 1,2,3,...$ this function is increasing in the parameter $m$ and decreasing in parameters $k$ and $\theta$. In Theorems 4-5 below we derive the cumulants and moments, framed in terms of this function for brevity. These theorems generalise cumulant and moment formulae shown for the classical coupon-collector distribution in Read (1998, p. 176).[3] Our formulae extend the results in previous literature to the more general case where $k$ and $\theta$ are free parameters.

**THEOREM 4 (Cumulants):** The Maclaurin expansion for the cumulant generating function is:

$$K(s) = -ks + \sum_{n=1}^{\infty} \frac{(1 - e^{-s})^n}{n} \cdot h_n(m, k, \theta),$$

with a radius of convergence that includes a neighbourhood of $s = 0$. Consequently, the $r$th cumulant of the negative occupancy distribution is:

$$\kappa_r = \sum_{i=1}^{r} (-1)^{r-i} \cdot S(r, i) \cdot (i - 1)! \cdot h_i(m, k, \theta) - k \cdot \mathbb{I}(r = 1).$$

**COROLLARY (Mean, variance, skewness, kurtosis):** The first four cumulants of the negative occupancy distribution are:

$$\mu_{m,k,\theta} \equiv \kappa_1 = h_1(m, k, \theta) - k,$$

$$\sigma_{m,k,\theta}^2 \equiv \kappa_2 = h_2(m, k, \theta) - h_1(m, k, \theta),$$

$$\kappa_3 = 2h_3(m, k, \theta) - 3h_2(m, k, \theta) + h_1(m, k, \theta),$$

$$\kappa_4 = 6h_4(m, k, \theta) - 12h_3(m, k, \theta) + 7h_2(m, k, \theta) - h_1(m, k, \theta).$$

The first two of these cumulants are the mean and variance of the distribution. The skewness and kurtosis of the distribution are obtained from the cumulants as:

$$\gamma_{m,k,\theta} \equiv \mathbb{S}\mathbb{k}\mathbb{e}\mathbb{w}(T_k) = \frac{\kappa_3}{\sigma_{m,k,\theta}^3} \qquad \kappa_{m,k,\theta} \equiv \mathbb{K}\mathbb{u}\mathbb{r}\mathbb{t}(T_k) = 3 + \frac{\kappa_4}{\sigma_{m,k,\theta}^4}.$$

---

[3] That paper examined the distribution of the hitting times, not the *excess* hitting times, so the cumulant function did not include the $-ks$ term, and the resulting expression for the mean did not include the location shift term. The remaining cumulants can be obtained from our formula by substituting $k = m$ and $\theta = 1$.

It is again important to stress that the negative occupancy distribution models the *excess* hitting time for the occupancy process, not the total hitting time. Since $T_k$ represents the number of excess balls to get to occupancy number $k$ (i.e., $k$ out of $m$ bins occupied), the total number of balls required for this occupancy number is $T_k + k$. This latter random variable has a location shift, giving it the mean and variance:

$$\mathbb{E}(T_k + k) = h_1(m, k, \theta),$$
$$\mathbb{V}(T_k + k) = h_2(m, k, \theta) - h_1(m, k, \theta).$$

In the classical coupon-collector problem we examine the case where $\theta = 1$ and we want the hitting time to full occupancy (i.e., $k = m$). In this case, we have the well-known formulae:

$$\mathbb{E}(T_m + m) = m H_m^{(1)},$$
$$\mathbb{V}(T_m + m) = m^2 H_m^{(2)} - m H_m^{(1)}.$$

Computation of the moments of the negative occupancy distribution requires computation of the generalised harmonic numbers. For large $m$ it is possible to approximate the moments by their limiting forms. By bounding the sum of powers in the generalised harmonic numbers by upper and lower integrals, and then applying the limit with the squeeze theorem, it can be shown that $h_n(m, k, \theta) \to k/\theta^n$ as $m \to \infty$. Thus, the cumulants converge to:

$$\lim_{m \to \infty} \kappa_r = k \left[ \sum_{i=1}^{r} (-1)^{r-i} \cdot S(r, i) \cdot (i-1)! \cdot \theta^{-i} - \mathbb{I}(r = 1) \right].$$

This gives the asymptotic moments:

$$\mu_{\infty, k, \theta} = k \cdot \frac{1 - \theta}{\theta} \qquad \sigma^2_{\infty, k, \theta} = k \cdot \frac{1 - \theta}{\theta^2},$$

$$\gamma_{\infty, k, \theta} = \frac{1}{\sqrt{k}} \cdot \frac{2 - \theta}{\sqrt{1 - \theta}} \qquad \kappa_{\infty, k, \theta} = 3 + \frac{1}{k} \cdot \left( 6 + \frac{\theta^2}{1 - \theta} \right).$$

These are the moments of the negative binomial distribution with probability parameter $1 - \theta$, which is unsurprising, since the negative occupancy distribution converges to this distribution in the limit. While this limiting approximation gives a crude approximation to the moments in the more general case, it can be improved by looking at more specific limits. In particular, it is useful to see what happens to the distribution and its moments as $m$ and $k$ both become large, but the ratio $k/m$ approaches a fixed limit. The simplest way to determine this is to derive the limiting form of the cumulant function, and use this to find the limiting values of the cumulants and central moments of interest. This is done below in Theorems 5-6, which yield a more generalised asymptotic form that approaches the above limits.

**THEOREM 5 (Asymptotic cumulant generating function):** If $m \to \infty$ and $k \to \infty$ with a fixed ratio $k/m \to \lambda$ then the cumulant function has the asymptotically equivalent limiting form:

$$K(s) \sim m \begin{bmatrix} \lambda \ln(\theta) - (1-\lambda) \ln|1-\lambda| \\ -\dfrac{1-(1-\theta)e^s}{\theta e^s} \cdot \ln|1-(1-\theta)e^s| \\ +\dfrac{1-(1-\theta)e^s - \lambda \theta e^s}{\theta e^s} \cdot \ln|1-(1-\theta)e^s - \lambda \theta e^s| \end{bmatrix}.$$

In the classical case where $\theta = 1$ we have:

$$K(s) \sim m \left[ -(1-\lambda) \ln|1-\lambda| + \frac{1-\lambda e^s}{e^s} \cdot \ln|1-\lambda e^s| \right].$$

(Note that the parameter $m$ appears on the right-hand-side only as a multiplicative term, which means that we can also obtain a corresponding limit statement by dividing both sides by $m$.)

**THEOREM 6 (Asymptotic moments):** If $m \to \infty$ and $k \to \infty$ with $k/m \to \lambda$ then we have the following asymptotically equivalent forms for the mean and variance of the distribution:

$$\mu_{m,k,\theta} \sim \mu^*_{m,k,\theta} = -k - \frac{m}{\theta} \cdot \ln\left(\frac{m-k}{m}\right),$$

$$\sigma^2_{m,k,\theta} \sim \sigma^{*2}_{m,k,\theta} = \frac{m}{\theta^2} \cdot \frac{k}{m-k} + \frac{m}{\theta} \cdot \ln\left(\frac{m-k}{m}\right).$$

All higher order cumulants are of the form $\kappa^*_r = m \cdot f(\lambda, \theta)$, which means that all the higher order cumulants of $T_k/m$ are of order $O(m^{-1})$, so they vanish in the limit. (We show the form of the asymptotic third and fourth moments in the proof, but omit them here.)

The above results show that the higher order cumulants of $T_k/m$ vanish in the limit, and so the distribution of this scaled excess hitting time converges to the normal distribution. Since this scaled variable has non-negative support, it can be fruitfully approximated by any continuous distribution with non-negative support that converges to the normal distribution. Read (1998) argues in favour of a lognormal approximation for the distribution of the total hitting time in the coupon-collector problem, but we find the gamma distribution to be a better approximation. We will develop this approximation more specifically in Section 3, and examine its accuracy a large matrix of values for the space parameter and occupancy parameter. As expected from Theorems 5-6, we find that the approximation becomes highly accurate as $m$ and $k$ get large. This will give us a relatively simple method of computing the distribution in cases where large parameter values make other methods too computationally intensive.

## 3. Computing the negative occupancy distribution

The probability mass function for the negative occupancy distribution involves the noncentral Stirling numbers of the second kind, which present some computational challenges. Attempts to compute directly from the mass function or the weighted geometric representation run into difficulties with arithmetic underflow, and cancelling of positive and negative terms. One can attempt to deal with this difficulty by splitting these forms into their positive and negative parts (i.e., splitting the sums over the index values that give positive and negative terms), computing each part separately in log-space, and then adding them to get the desired value. Unfortunately, both the positive and negative parts of the sum tend to be very large relative to their difference, and so even very small amounts of rounding error due to arithmetic underflow can lead to a situation in which the negative part exceeds the positive part, giving a negative value for the mass function. (As a test, the present author has attempted some variations on computing the mass function by this method, and was unsuccessful in obtaining a stable algorithm that does not give negative values of the mass function for some arguments.)

Rather than attempting direct computation from the formulae for these numbers, it is preferable to compute the mass function recursively for occupancy values $r = 1, \ldots, k$ up to the desired occupancy number. We have already seen that with unit occupancy parameter the negative occupancy distribution reduces to the geometric distribution, and so it can be easily computed. In the theorem below we present a recursive formula to build up from this distribution to get to the desired occupancy number. To avoid arithmetic underflow we undertake all computations in log-space, so we give the recursive equation in log-space.

**THEOREM 7 (Recursive computation):** Let $L(t, k) \equiv \log \text{NegOcc}(t|m, k, \theta)$ denote the log-probabilities from the negative occupancy distribution, and define the quantities:

$$L_k \equiv \log\left(1 - \theta \cdot \frac{m-k}{m}\right).$$

The log-probability values satisfy the following recursive equation:

$$L(t, 1) = \log(\theta) + t \log(1 - \theta),$$

$$L(t, k+1) = \log\left(\theta \cdot \frac{m-k}{m}\right) + \text{logsumexp}\begin{pmatrix} L(t, k), \\ L_k + L(t-1, k), \\ 2 \cdot L_k + L(t-2, k), \\ \vdots \\ t \cdot L_k + L(0, k) \end{pmatrix}.$$

We can use the recursive equation in the above theorem to form an algorithm to generate the log-probabilities from the negative occupancy distribution. This is done in Algorithm 1 below in pseudo-code. Since the support of the distribution is unbounded, our algorithm ranges over all argument values up to a stipulated maximum value $T$. The algorithm fixes the values $m$ and $\theta$ and begins by forming a $(T+1) \times k$ matrix of values, with the rows giving the argument values $t = 0, ..., T$ and the columns giving the occupancy parameter values $r = 1, ..., k$. The algorithm computes the first row of the matrix by using the log-probabilities for the geometric distribution, and then computes each subsequent row using the recursive equation in the above theorem. The final row of the matrix gives the log-probabilities for the negative occupancy distribution with the occupancy parameter $k$. The algorithm returns either the log-probabilities or probabilities of the negative occupancy distribution.

```
                    ALGORITHM 1: Negative Occupancy Distribution

Input:      Maximum argument T (non-negative integer)
            Space parameter m (positive integer or Inf)
            Occupancy parameter k (positive integer no greater than m)
            Probability parameter θ (probability value)
            Logical value log (specifying whether output is log-probability)
Output:     Vector of probabilities/log-probabilities
            from the negative occupancy mass function for t = 0,…,T
```

```
#Deal with special case where m = Inf
if (m = Inf) {
  LNEGOCC <- log(NegBin(0:T|k,1-θ) }

#Deal with special case where θ = 0 or k = 0
if ((θ = 0)|(k = 0)) {
  LNEGOCC <- [0, -Inf, … , -Inf] (with length T+1) }

#Deal with remaining cases
if ((θ > 0)&(m < Inf)&(0 < k ≤ m)) {
  #Generate empty matrix of log-probabilities
  L   <- Matrix with rows 0,…,T and columns 1,…,k
         (Initial values all set to -∞)
  LLL <- (1:k)*log(1 - θ*(m-(1:k))/m)

  #Compute the first row
  for each t = 1,…,T {
    L[t,1] <- log(θ) + t*log(1-θ)
    for each r = 2,…,k {
      L[t,r] <- log(θ*(m-r+1)/m) +
        logsumexp((0:t)*LLL[k] + L[t:0,r-1]) } }
  LNEGOCC <- L[, k] }

#Return output
if (log) { LNEGOCC } else { exp(LNEGOCC) }
```

Algorithm 1 allows the user to compute the mass function for a given parameter specification for all arguments up to a stipulated maximum. The computational intensity is proportionate to the size of the occupancy parameter $k$. In the course of computing the mass function over the argument values $t = 0, \ldots, T$ the algorithm computes the mass function for all values of the occupancy parameter over $r = 1, \ldots, k$. The downside of the algorithm is that the recursive method means that these other intermediate computations are generally wasted. However, if we store these intermediate values instead of discarding them, then the method allows the user to compute "blocks" of negative occupancy mass values for all occupancy parameters up to the space parameter. To illustrate the accuracy of the algorithm, we compute the negative occupancy distribution with parameters $m = 30$, $k = 14$ and $\theta = 0.6$ and simulate from the distribution using sums of geometric random variables using the convolution representation. A bar-plot of $10^7$ simulations and exact probabilities from the mass function is shown in Figure 1 below. This confirms the accuracy of the algorithm, against large-sample simulation.

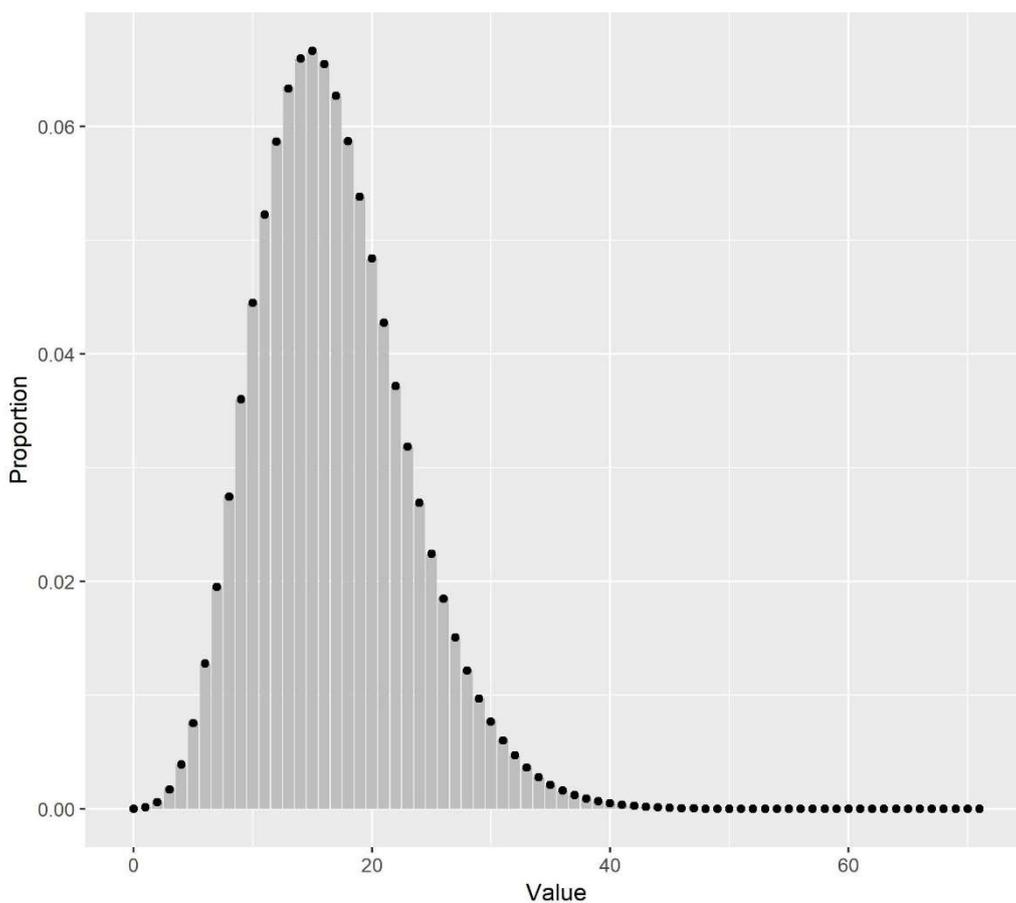

**FIGURE 1:** Bar-plot of simulations from the negative occupancy distribution

This figure uses $10^7$ simulations with parameters $m = 30$, $k = 14$ and $\theta = 0.6$.
Simulated proportions are the grey bars; the true mass values are the black dots.

We will examine the accuracy of the gamma distribution as an approximation to the negative occupancy distribution. Using a continuous distribution to approximate a discrete distribution opens up various specific ways in which the approximation may be formulated. In our analysis we use an approximation that treats each discrete argument value in the negative occupancy distribution as an equal-length interval of the continuous approximating value. Using the fact that $T_k/m$ converges to a gamma distribution (which converges to a normal distribution) under the limit conditions in Theorem 6, we set the approximation as:

$$\widehat{\text{NegOcc}}(t|m,k,\theta) \equiv \int_t^{t+1} \text{Ga}(r|\alpha_{m,k,\theta}, \beta_{m,k,\theta})\,dr,$$

where the parameters of the gamma distribution are given by:[4]

$$\alpha_{m,k,\theta} = \frac{(\mu_{m,k,\theta} + \tfrac{1}{2})^2}{\sigma^2_{m,k,\theta}} \qquad \beta_{m,k,\theta} = \frac{\mu_{m,k,\theta} + \tfrac{1}{2}}{\sigma^2_{m,k,\theta}}.$$

To allow accurate computation for small probabilities, and avoid underflow problems, we will conduct our computation in log-space. For brevity, let $L_{\text{Ga}}$ denote the log-distribution function for the gamma distribution with the above parameters (i.e., giving logarithms of the cumulative probabilities from this distribution). We form our approximation as:[5]

$$\hat{L}(t,k) = \text{logdiffexp}\begin{pmatrix} L_{\text{Ga}}(t+1), \\ L_{\text{Ga}}(t) \end{pmatrix} \qquad L_{\text{Ga}}(t) \equiv \log \int_0^t \text{Ga}(r|\alpha_{m,k,\theta}, \beta_{m,k,\theta})\,dr.$$

---

[4] This parameterisation of the gamma distribution includes a "continuity correction" for the transition from a continuous approximating random variable to a discrete random variable. To derive the parameterisation, we let $S \sim \text{Ga}(\alpha, \beta)$ and take $T_k = \lfloor S \rfloor$ as our approximation. Using the continuity correction $\lfloor S \rfloor \approx S - \tfrac{1}{2}$ gives:

$$\mu = \mathbb{E}(T_k) = \mathbb{E}(\lfloor S \rfloor) \approx \mathbb{E}(S - \tfrac{1}{2}) = \frac{\alpha}{\beta} - \tfrac{1}{2},$$

$$\sigma^2 = \mathbb{V}(T_k) = \mathbb{V}(\lfloor S \rfloor) \approx \mathbb{V}(S - \tfrac{1}{2}) = \frac{\alpha}{\beta^2}.$$

This yields (approximating) simultaneous equations relating the moments of the negative occupancy distribution to the parameters of the (approximating) gamma distribution:

$$\mu + \tfrac{1}{2} = \frac{\alpha}{\beta} \qquad \sigma^2 = \frac{\alpha}{\beta^2}.$$

Solving for the parameters of the gamma distribution gives:

$$\alpha = \frac{(\mu + \tfrac{1}{2})^2}{\sigma^2} \qquad \beta = \frac{\mu + \tfrac{1}{2}}{\sigma^2}.$$

The parameterisation shown in the body of the paper uses these formulae, but included subscripts in the notation to show the dependence on the underlying parameters of the negative occupancy distribution.

[5] The function logdiffexp is defined by $\text{logdiffexp}(\ell_1, \ell_2) \equiv \log(\exp(\ell_1) - \exp(\ell_2))$; it computes differences in log-space. It is analogous to the logsumexp function, but it computes differences instead of sums.

We can use the above equation to form an algorithm to generate approximate log-probabilities from the negative occupancy distribution. This is done in Algorithm 2 below in pseudo-code. As with our previous algorithm for computing the exact mass function, since the support of the distribution is unbounded, our algorithm ranges over all argument values up to a stipulated maximum value $T$. The algorithm computes the mean and variance of the negative occupancy distribution and uses these to compute the parameters of the gamma approximation. The approximate log-probabilities are then computed in a loop that preserves the previous value of the log-distribution function to improve computational efficiency.

```
          ALGORITHM 2: Negative Occupancy Distribution (Approx)

Input:       Positive integer T
             Positive integers m and k
             Real value θ
             Logical value log
Output:      Vector of negative occupancy log-probabilities for t = 0,…,T
```

```
#Create output vector
APPROX <- Vector with elements t = 0,…,T
          (Initial values all set to -Inf)

#Compute parameters
H1    <- sum(1/((m-k+1):m))
H2    <- sum(1/((m-k+1):m)^2)
MEAN  <- (m/θ)*H1 - k
VAR   <- (m/θ)^2*H2 - (m/θ)*H1
SHAPE <- (MEAN+1/2)^2/VAR
RATE  <- (MEAN+1/2)/VAR

#Approximation using gamma distribution
#The function logCDFGa is the gamma log-CDF
if (VAR = 0) { APPROX[0] <- 0 }
if (VAR > 0) {
  t     <- 0
  UPPER <- logGaF(0, rate = RATE, shape = SHAPE)
  while (t <= T) {
    LOWER <- UPPER
    UPPER <- logCDFGa(t+1, rate = RATE, shape = SHAPE)
    APPROX[t] <- logdiffexp(UPPER, LOWER)
    t <- t+1 }

#Return output
if (log) { APPROX } else { exp(APPROX) }
```

Using the same parameterisation for the negative occupancy distribution as in Figure 1 above, we can compute the gamma approximation to this distribution using Algorithm 2 and compare this to the exact mass values computed from Algorithm 1. A plot of the exact probabilities and the approximation is shown in Figure 2 below.

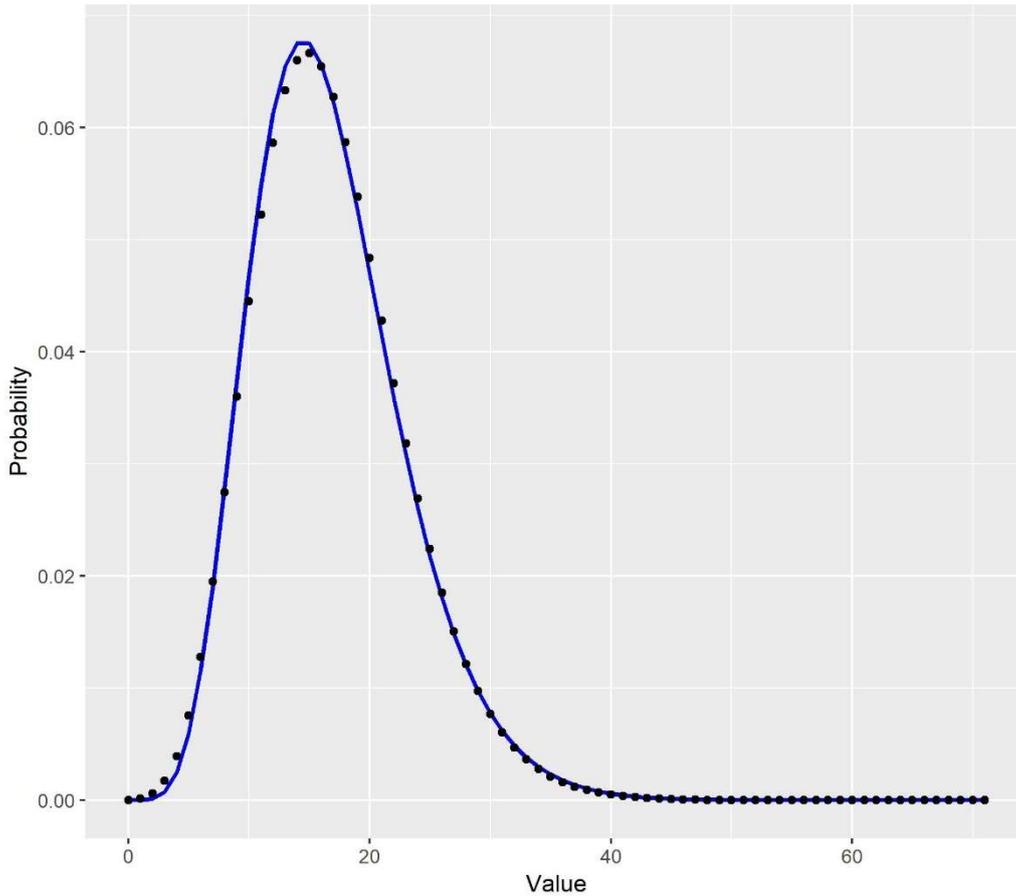

**FIGURE 2:** Gamma approximation to the negative occupancy distribution

Negative occupancy distribution has parameters $m = 30$, $k = 14$ and $\theta = 0.6$.
Gamma approximation is the blue line; the true mass values are the black dots.

We can see from the figure that the gamma distribution provides a reasonable approximation to the negative occupancy distribution even in the present case, where the parameter values are small. In this case, the approximation has greater positive skewness than the true distribution, so it underestimates the true probabilities in the left-tail and overestimates the true probabilities in the right-tail. Values in the "main body" of the distribution are overestimated to the left and underestimated to the right. (To the naked eye the approximation looks to be accurate in the tails of the distribution, but this is just an artefact of the small values in the tails. Plotting the exact and approximate log-probabilities shows that the approximation overestimates the probabilities in the right-tail, though both sets of values are extremely small in any case.) It is simple to create a hybrid algorithm that either uses recursive computation (Algorithm 1) or the gamma approximation (Algorithm 2) depending on the preference of the user. In general, the approximation will be okay when the space parameter $m$ is "large" and it becomes necessary when this parameter is "too large" to make recursive computation feasible.

# 4. Accuracy of the approximation

In order to obtain a more systematic idea of the relative approximation accuracy as we increase the parameters, we will measure the accuracy of our approximation over a large block of values of the space and occupancy parameters. Since Algorithm 1 computes the negative occupancy distribution recursively over the occupancy parameter, it is already efficient for computing the log-probabilities of the negative occupancy distribution over a block of parameter values. For simplicity, we examine the classical case where $\theta = 1$ and we compute the negative occupancy distribution and its approximation over a block of parameter values $0 < k \leq m \leq M$.

Our analysis here is similar to a corresponding analysis in O'Neill (2020), which computed the accuracy of the normal approximation to the classical occupancy distribution, taken over a large block of its size and space parameters. That paper used the log-root-mean-square-error between the true distribution and its approximation to assess the accuracy of the approximation. In the present context that measure of accuracy is invalid, since an average error taken over the infinite support of the negative occupancy distribution will be zero. In the present context we will instead measure the accuracy using the root-square-error:

$$\text{RSE}(m, k) \equiv \left( \sum_{t=0}^{\infty} (\text{NegOcc}(t|m, k, 1) - \widehat{\text{NegOcc}}(t|m, k, 1))^2 \right)^{1/2}.$$

The root-squared error sums over the entire support of the distribution, so the upper bound in the sum is infinity. Unfortunately, it is not possible to compute the true distribution or its approximation over all the values in the support, so we instead (under) estimate the accuracy using a truncated sum that we cut-off at a large finite value $T(m, k)$, which is computed based on the parameters of the distribution.[6] Truncation reduces the value of the accuracy measure, but since the values in the sum vanish rapidly when both distributions give small probability values, it still gives an accurate approximation to the true distance.

$$\widehat{\text{RSE}}(m, k) \equiv \left( \sum_{t=0}^{T(m,k)} (\text{NegOcc}(t|m, k, 1) - \widehat{\text{NegOcc}}(t|m, k, 1))^2 \right)^{1/2}.$$

---

[6] For our analysis of the LRSE we choose the cut-off to be the smallest integer that is at least fie standard deviations above the mean —i.e., we take $T(m,k) = \lceil \mu_{m,k,1} + 5\sigma_{m,k,1} \rceil$ for all $0 < k \leq m \leq M$. For each $m = 1, \ldots, M$ we use a variation of Algorithm 1 to compute the block of log-probabilities for $k = 1, \ldots, m$ with arguments taken over the values $t = 0, \ldots, \max_k T(m, k)$. We compute each value $\widehat{\text{RSE}}$ using $t = 0, \ldots, T(m, k)$, so some of the log-probabilities in our block of values remain unused. This method is used for computational efficiency, but still ensuring consistency of the accuracy measure over the different parameter values in the block.

The RSE measures the vector distance between probability vectors for the negative occupancy distribution and its approximation. Letting $\mathbf{p}_{m,k,\theta} = (p_{t,m,k,\theta}|t=0,1,2,...)$ denote the vector of negative occupancy mass values and $\hat{\mathbf{p}}_{m,k,\theta} = (\hat{p}_{t,m,k,\theta}|t=0,1,2,...)$ denote the vector of gamma mass values, we have $\text{RSE}(m,k) = \|\hat{\mathbf{p}}_{m,k,\theta} - \mathbf{p}_{m,k,\theta}\|$. Our estimate $\widehat{\text{RSE}}$ likewise measures the vector distance of the truncated vectors.

Figure 3 shows a heatmap of the root-squared-error (RSE) of the gamma approximation to the negative occupancy distribution, taken over a block of parameter values up to $M = 1000$. High accuracy is represented by the darker areas showing low values of the RSE. The plot shows that there are ridges of high and low accuracy, corresponding to particular values of the ratio $k/m$. The accuracy of the gamma approximation also generally increases (i.e., lower RSE) as the space and occupancy parameters become large.

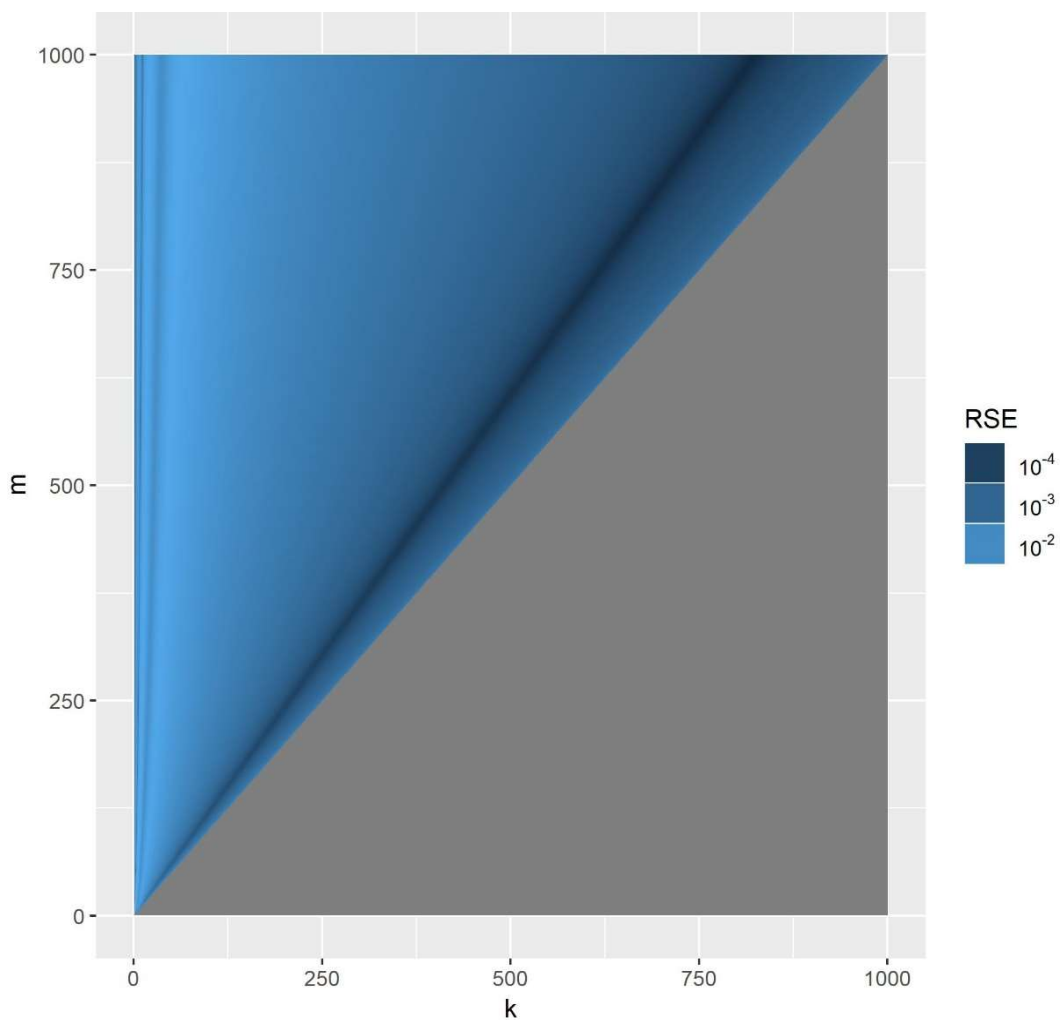

**FIGURE 3:** Heatmap of RSE values for gamma approximation

Negative occupancy distribution with parameters $0 < k \leq m \leq 1000$ and $\theta = 1$.

The heatmap in Figure 3 gives a general overview of the accuracy, taken over all combinations of space and occupancy parameters in a large block. It is useful to reduce this to one-dimension by looking at the maximum RSE, mean RSE and diagonal RSE for each value of $m$ in the block of parameters. This is shown in Figure 4 below. We can see that the mean accuracy increases (i.e., mean RSE decreases) as we increase the value of $m$, but the accuracy of the least accurate approximation (i.e., the maximum RSE) appears to plateau. This occurs because the RSE map exhibits a "ridge" of high values for low $k$, where the accuracy of the gamma approximation does not improve with higher $m$. The diagonal values in the heatmap are shown in the bottom line of Figure 4, giving the accuracy of the approximation to the coupon-collector distribution. We can see that both the mean and diagonal RSE appear to decrease slower than exponential decay. (Since the RSE is shown on a log-scale, exponential decay would be a straight line.) Using a space parameter of $m \geq 289$ is sufficient to get mean RSE less than one-percent.

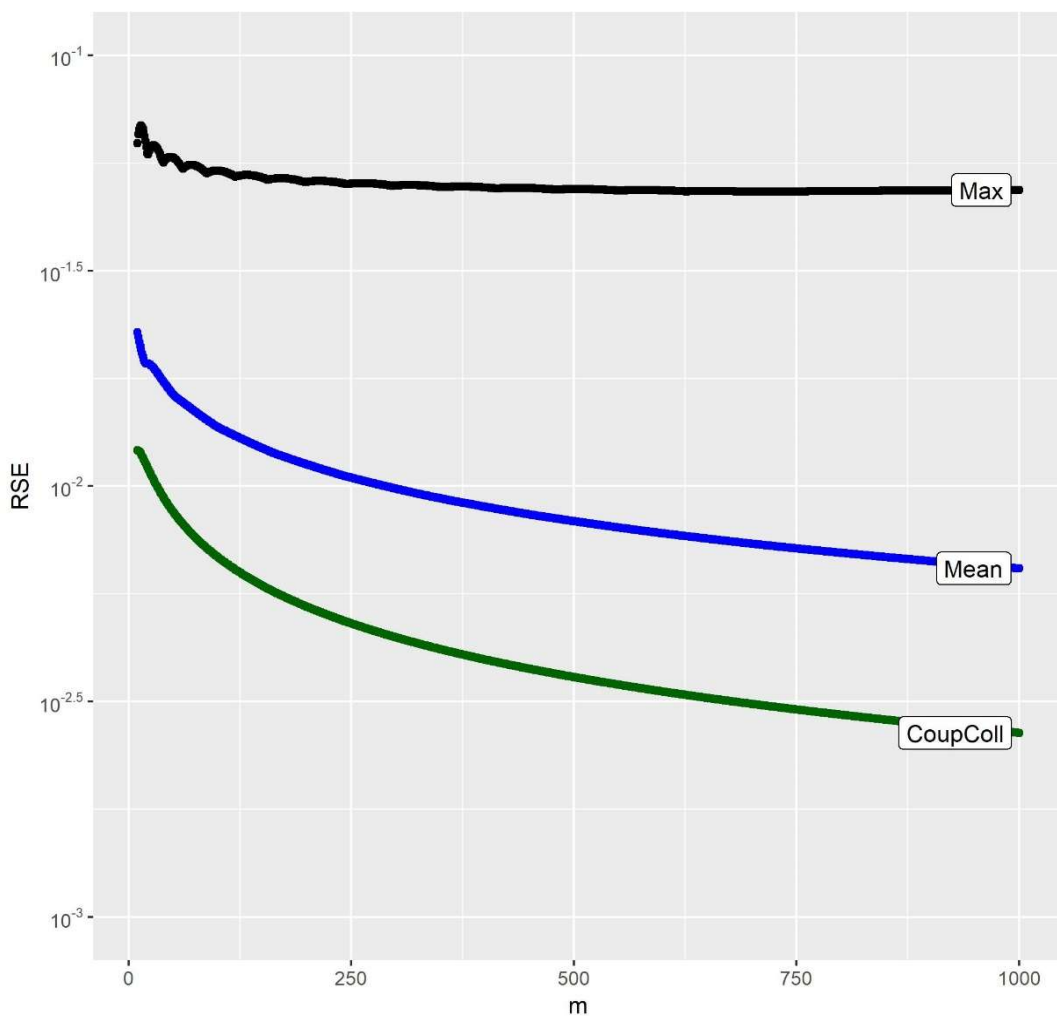

**FIGURE 4:** Plot of RSE values for gamma approximation

Negative occupancy distribution with parameters $10 \leq m \leq 1000$ and $\theta = 1$.

Algorithm 1 gives a general method for computing the mass function of the negative occupancy distribution, that is limited only by the computational power and storage space of the machine. However, as $m$ increases, the computational burden of the recursive algorithm increases, and at the same time, the accuracy of the gamma approximation gets better and better, making the latter more desirable. Thus, to develop a reliable and computationally efficient algorithm for the mass function of the distribution, the simplest method is to use exact recursive computation in Algorithm 1 for all values of $m$ that are "not too large", and then switch to the approximation formula when the value of $m$ becomes "too large".

It remains to decide when one should make the switch from the exact recursive method to the gamma approximation. Since this is a trade-off between accuracy and computational burden, there is no sacrosanct answer, and it is a matter of deciding how much computation time and storage space is available, and how much accuracy is required. The methods in this paper are implemented in the `occupancy` package in `R` (O'Neill 2021b), which has automated functions to compute the probability mass function, cumulative distribution function, quantile function, and simulation function for the negative occupancy distribution. The probability functions in this package allow "exact" computation by the recursive method, or approximate computation using the gamma approximation shown in this paper.

## 5. Conclusion

In this paper we have derived a number of properties of the negative occupancy distribution and coupon-collector distributions introduced in O'Neill (2021a). These distributions arise to describe the behaviour of the excess hitting times in the extended occupancy problem. They can be used to solve the classical "coupon collector's problem", and extensions of this problem involving cases where the collector seeks only a partial collection, or cases where coupons may be damaged, etc. The family of negative occupancy distributions is indexed by parameters $m$, $k$ and $\theta$, taken from the specification of the extended occupancy problem.

In the present paper we have showed that the negative occupancy distribution is obtains as the convolution of geometric distributions with a specified set of parameters. We have derived the generating functions and moments of the distributions, and asymptotic forms for these in the case where the parameters $m$ and $k$ are large. We have also advanced the view that the gamma distribution provides a good continuous approximation to the distribution in this case, and we have examined the accuracy of this distribution over a large block of parameter values. To do this, we have developed a recursive algorithm to compute the mass function of the negative occupancy distribution, and an algorithm to compute its approximation.

The negative occupancy distribution is a generalisation of the negative binomial distribution, with an additional parameter the limits the number of "bins" in the occupancy problem to a finite value. Using a finite number of bins means that the balls in the occupancy problem can be allocated to bins that are already occupied, which means that the number of occupied bins tends to be lower, and the excess hitting time for any given number of bins tends to be higher. The distribution arises directly in marginal and conditional occupancy problems, but it also has potential uses more generally in count models, particularly in cases where the user desires an extension to the negative binomial distribution.

We hope that this paper has piqued the interest of the reader in the extended occupancy problem and the interesting families of discrete distributions that arise from this problem. The negative occupancy distribution under analysis in this paper is an interesting family of distributions, and it has a number of interesting mathematical connections to the Stirling numbers of the second kind and the generalised harmonic numbers. Functions for this probability distribution are available in the `occupancy` package in `R`, so it is simple for practitioners to implement this distribution in their work if desired.

**Declarations**

| | |
|---|---|
| **Availability of data and materials** | Not applicable |
| **Competing interests** | The author declares that he has no competing interests |
| **Funding** | There was no specific funding for this paper or project |
| **Author's contributions** | The author undertook all work on the present paper (i.e., conception, analysis and writing for the paper). |
| **Acknowledgements** | The author would like to thank the *Journal for Statistical Distributions and Applications* for review of the paper. |

# Appendix: Proof of Theorems

**PROOF OF THEOREM 1:** For the base case $k = 1$ the noncentral Stirling numbers have the form $S(n, 0, \phi) = \phi^n$ which gives the mass function:

$$\text{NegOcc}(t|m, 1, \theta) = \frac{\theta^{t+1}}{m^{t+1}} \cdot m \cdot S\left(t, 0, m \cdot \frac{1-\theta}{\theta}\right)$$

$$= \frac{\theta^{t+1}}{m^{t+1}} \cdot m \cdot \left(m \cdot \frac{1-\theta}{\theta}\right)^t$$

$$= (1-\theta)^t \theta = \text{Geom}(t|\theta).$$

To establish the recursive equation, we will use the "telescoping sum" rule for the noncentral Stirling numbers (O'Neill 2020a, Appendix I, Lemma 2), which gives:

$$S\left(k + t, k, m \cdot \frac{1-\theta}{\theta}\right) = \sum_{r=0}^{t} S\left(k + t - r - 1, k - 1, m \cdot \frac{1-\theta}{\theta}\right) \cdot \left(k + m \cdot \frac{1-\theta}{\theta}\right)^r.$$

Using this telescoping sum we then have:

$$\text{NegOcc}(t|m, k+1, \theta) = \frac{\theta^{k+t+1}}{m^{k+t+1}} \cdot (m)_{k+1} \cdot S\left(k + t, k, m \cdot \frac{1-\theta}{\theta}\right)$$

$$= \frac{\theta^{k+t+1}}{m^{k+t+1}} \cdot (m)_{k+1} \sum_{r=0}^{t} S\left(k + t - r - 1, k - 1, m \cdot \frac{1-\theta}{\theta}\right) \cdot \left(k + m \cdot \frac{1-\theta}{\theta}\right)^r$$

$$= \frac{\theta^{k+t+1}}{m^{k+t+1}} \cdot (m)_{k+1} \sum_{s=0}^{t} S\left(k + s - 1, k - 1, m \cdot \frac{1-\theta}{\theta}\right) \cdot \left(k + m \cdot \frac{1-\theta}{\theta}\right)^{t-s}$$

$$= \sum_{s=0}^{t} \left[ \begin{array}{l} \frac{\theta^{k+s}}{m^{k+s}} \cdot (m)_k \cdot S\left(k + s - 1, k - 1, m \cdot \frac{1-\theta}{\theta}\right) \\ \times \frac{\theta^{t-s+1}}{m^{t-s+1}} \cdot (m - k) \cdot \left(k + m \cdot \frac{1-\theta}{\theta}\right)^{t-s} \end{array} \right]$$

$$= \sum_{s=0}^{t} \left[ \begin{array}{l} \frac{\theta^{k+s}}{m^{k+s}} \cdot (m)_k \cdot S\left(k + s - 1, k - 1, m \cdot \frac{1-\theta}{\theta}\right) \\ \times \left(1 - \theta \cdot \frac{m-k}{m}\right)^{t-s} \cdot \left(\theta \cdot \frac{m-k}{m}\right) \end{array} \right]$$

$$= \sum_{s=0}^{t} \text{NegOcc}(s|m, k, \theta) \cdot \text{Geom}\left(t - s \middle| \theta \cdot \frac{m-k}{m}\right)$$

$$= \left[\text{NegOcc}(m, k, \theta) * \text{Geom}\left(\theta \cdot \frac{m-k}{m}\right)\right](t).$$

Since this recursive result holds for all arguments $t$, this establishes the recursive convolution equation in the theorem. ∎

**PROOF OF THEOREM 2:** This result uses a variant of the weighted geometric representation in Sen and Balakrishnan (1999, Theorem 1, p. 422). That theorem looked in general terms at sums of geometric random variables on the support 1,2,3, ..., but we will adjust it here to give the corresponding result for sums of geometric random variables over the support 0,1,2,3, .... Making this adjustment, if we have a sequence of geometric random variables $S_1, S_2, S_3, ...$ with $S_i \sim \text{Geom}(q_i)$ (in our notation) then their sum $T_k = \sum_{\ell=1}^{k} S_k$ has probability mass function:

$$\mathbb{P}(T_k = t) = \sum_{i=1}^{k} w_{i,k} \cdot \text{Geom}(t + k - 1 | q_i) \qquad w_{i,k} = \prod_{j=1, j \neq i}^{k} \frac{q_j}{q_j - q_i}.$$

This result requires the values $q_1, q_2, q_3, ...$ to be distinct, which is satisfied in our case. In our case we have $q_i = \theta(m - k + 1)/m$, which establishes the main equation of the theorem and gives the weights:

$$w_{i,k} = \prod_{j=1, j \neq i}^{k} \frac{q_j}{q_j - q_i} = \prod_{j=1, j \neq i}^{k} \frac{(m - j + 1)}{(m - j + 1) - (m - i + 1)}$$

$$= \prod_{j=1, j \neq i}^{k} \frac{(m - j + 1)}{i - j}$$

$$= \frac{(m)_k}{(m - i + 1) \prod_{j=1, j \neq i}^{k} (i - j)}.$$

Splitting the product term in the denominator (with index implicitly over $1 \leq j \leq k$) into its positive and negative parts, gives:

$$\prod_{j \neq i} (i - j) = \prod_{j < i} (i - j) \prod_{j > i} (i - j)$$

$$= (-1)^{k-i} \prod_{j < i} (i - j) \prod_{j > i} (j - i)$$

$$= (-1)^{k-i} (i - 1)! \, (k - i)!$$

We therefore have:

$$w_{i,k} = \frac{(m)_k}{(m - i + 1) \prod_{j=1, j \neq i}^{k} (i - j)} = (-1)^{k-i} \frac{(m)_k}{(m - i + 1)(i - 1)! \, (k - i)!}.$$

This establishes the weighted geometric form and the corresponding weight values given in the theorem, which is enough to complete the proof. However, for the sake of completeness, we will also show that this gets us back to the standard form for the mass function of the negative occupancy distribution. (This can also be regarded as an alternative proof.) Using the weighted geometric representation asserted in the theorem have:

$$\text{NegOcc}(t|m,k,\theta) = \sum_{i=1}^{k} w_{i,k} \cdot \left(1 - \theta \cdot \frac{m-i+1}{m}\right)^{t+k-1} \left(\theta \cdot \frac{m-i+1}{m}\right)$$

$$= \sum_{i=1}^{k} w_{i,k} \cdot \left(\frac{\theta}{m} \cdot \left(\frac{m}{\theta} - (m-i+1)\right)\right)^{t+k-1} \left(\frac{\theta}{m} \cdot (m-i+1)\right)$$

$$= \sum_{i=1}^{k} w_{i,k} \cdot \left(\frac{\theta}{m} \cdot \left(i - 1 + m \cdot \frac{1-\theta}{\theta}\right)\right)^{t+k-1} \left(\frac{\theta}{m} \cdot (m-i+1)\right)$$

$$= \frac{\theta^{k+t}}{m^{k+t}} \cdot \sum_{i=1}^{k} w_{i,k} \cdot \left(i - 1 + m \cdot \frac{1-\theta}{\theta}\right)^{t+k-1} (m-i+1)$$

$$= \frac{\theta^{k+t}}{m^{k+t}} \cdot (m)_k \cdot \sum_{i=1}^{k} \frac{1}{(i-1)!\,(k-i)!} (-1)^{k-i} \left(i - 1 + m \cdot \frac{1-\theta}{\theta}\right)^{t+k-1}$$

$$= \frac{\theta^{k+t}}{m^{k+t}} \cdot (m)_k \cdot \frac{1}{(k-1)!} \sum_{i=1}^{k} \binom{k-1}{i-1} (-1)^{k-i} \left(i - 1 + m \cdot \frac{1-\theta}{\theta}\right)^{k+t-1}$$

$$= \frac{\theta^{k+t}}{m^{k+t}} \cdot (m)_k \cdot \frac{1}{(k-1)!} \sum_{i=0}^{k-1} \binom{k-1}{i} (-1)^{k-i-1} \left(i + m \cdot \frac{1-\theta}{\theta}\right)^{k+t-1}$$

$$= \frac{\theta^{k+t}}{m^{k+t}} \cdot (m)_k \cdot S\left(k+t-1, k-1, m \cdot \frac{1-\theta}{\theta}\right).$$

We therefore see that the representation can be obtained directly by algebraic manipulation of the mass function of the negative occupancy distribution. (Taking these algebraic steps in reverse order gives an alternative proof of the theorem that does not rely on the result in Sen and Balakrishnan.) ∎

**PROOF OF THEOREM 3:** We will prove the form of the probability generating function; all the other generating functions are simple variations obtained by substitution of argument values. To derive the probability generating function we use the ordinary generating function for the noncentral Stirling numbers of the second kind (Charalambides 2005, p. 93, Eqn. 2.40):

$$\text{OGF}(x|k,\phi) \equiv \sum_{n=k}^{\infty} S(n,k,\phi) \cdot x^n = \frac{x^k}{\prod_{\ell=0}^{k}(1 - (\ell + \phi)x)} \qquad |x| < \frac{1}{\phi + k}.$$

Using this result we obtain:

$$G(z) = \sum_{t=0}^{\infty} z^t \cdot \text{NegOcc}(t|m,k,\theta)$$

$$= \sum_{t=0}^{\infty} z^t \cdot \frac{\theta^{k+t}}{m^{k+t}} \cdot (m)_k \cdot S\left(k+t-1, k-1, m \cdot \frac{1-\theta}{\theta}\right)$$

$$= \frac{\theta}{m} \cdot \frac{(m)_k}{z^{k-1}} \cdot \sum_{t=0}^{\infty} \left(\frac{z\theta}{m}\right)^{k+t-1} \cdot S\left(k+t-1, k-1, m \cdot \frac{1-\theta}{\theta}\right)$$

$$= \frac{\theta}{m} \cdot \frac{(m)_k}{z^{k-1}} \cdot \text{OGF}\left(\frac{z\theta}{m} \bigg| k-1, m \cdot \frac{1-\theta}{\theta}\right)$$

$$= \frac{\theta}{m} \cdot \frac{(m)_k}{z^{k-1}} \cdot \frac{\left(\frac{z\theta}{m}\right)^{k-1}}{\prod_{\ell=0}^{k-1}\left(1 - \left(\ell + m \cdot \frac{1-\theta}{\theta}\right)\frac{z\theta}{m}\right)}$$

$$= \theta^k \cdot \frac{(m)_k}{m^k} \cdot \prod_{\ell=0}^{k-1} \frac{1}{1 - \left(\frac{m}{\theta} - (m-\ell)\right)\frac{z\theta}{m}}$$

$$= \theta^k \cdot \frac{(m)_k}{m^k} \cdot \prod_{\ell=0}^{k-1} \frac{m}{m - (m - (m-\ell)\theta)z}$$

$$= \theta^k \cdot \prod_{\ell=0}^{k-1} \frac{m-\ell}{m - (m - (m-\ell)\theta)z}$$

$$= \theta^k \cdot \prod_{\ell=m-k+1}^{m} \frac{\ell}{m - (m - \ell\theta)z}.$$

The domain of convergence for the ordinary generating function is:

$$\left|\frac{z\theta}{m}\right| < \frac{\theta}{m(1-\theta) + (k-1)\theta}.$$

Some simple algebra yields the equivalent requirement $|z| < m/(m - (m-k+1)\theta)$ and it can easily be shown that this inequality is the smallest upper bound that ensures that each of the inverse-product terms in the above expression is greater than zero. ∎

**ALTERNATIVE PROOF OF THEOREM 3:** For this proof we will use the fact that the negative occupancy distribution is a convolution of geometric distributions. Alternatively, if the reader is unwilling to accept this premise, the present proof can be regarded as a proof of that fact, by showing that the probability generating function of this convolution is the same as that of the negative occupancy distribution. In any case, observe that $T_k = \sum_{\ell=0}^{k-1} S_\ell$ and the increments in this sum each have probability generating function given by:

$$G_{S_\ell}(z) = \sum_{t=0}^{\infty} z^s \cdot \text{Geom}\left(s \bigg| \theta \cdot \frac{m-\ell}{m}\right)$$

$$\begin{aligned}
&= \sum_{t=0}^{\infty} z^s \cdot \left(1 - \theta \cdot \frac{m-\ell}{m}\right)^s \left(\theta \cdot \frac{m-\ell}{m}\right) \\
&= \left(\theta \cdot \frac{m-\ell}{m}\right) \sum_{t=0}^{\infty} \left(\left(1 - \theta \cdot \frac{m-\ell}{m}\right) \cdot z\right)^s \\
&= \frac{\theta \cdot \frac{m-\ell}{m}}{1 - \left(1 - \theta \cdot \frac{m-\ell}{m}\right) \cdot z} \\
&= \theta \cdot \frac{m-\ell}{m - (m - (m-\ell)\theta)z},
\end{aligned}$$

and the domain of convergence (required for convergence of the infinite geometric sum) is $|z| < m/(m - (m-\ell) \cdot \theta)$. We therefore have:

$$G(z) = \prod_{\ell=0}^{k-1} G_{S_\ell}(z) = \theta^k \cdot \prod_{\ell=0}^{k-1} \frac{m-\ell}{m - (m - (m-\ell)\theta)z} = \theta^k \cdot \prod_{\ell=m-k+1}^{m} \frac{\ell}{m - (m - \ell\theta)z},$$

with domain of convergence $|z| < m/(m - (m-k+1)\theta)$. This gives the asserted form for the probability generating function. The other generating functions are easily obtained from this result. ∎

**REMARK:** It is worth noting here that O'Neill (2021a, Theorem 7) shows that the limit $m \to \infty$ leads to $\text{NegOcc}(t|m, k, \theta) \to \text{NegBin}(t|k, 1-\theta)$. Thus, as a check on our working, we can confirm that our probability generating function leads to the probability generating function for the negative binomial distribution in the limit. Sure enough, we have:

$$\begin{aligned}
\lim_{m\to\infty} G(z) &= \theta^k \cdot \prod_{\ell=0}^{k-1} \lim_{m\to\infty} \frac{m-\ell}{m - (m - (m-\ell)\theta)z} \\
&= \theta^k \cdot \prod_{\ell=0}^{k-1} \lim_{m\to\infty} \frac{1 - \ell/m}{1 - (1 - (1-\ell/m)\theta)z} \\
&= \theta^k \cdot \prod_{\ell=0}^{k-1} \frac{1}{1 - (1-\theta)z} \\
&= \left(\frac{\theta}{1 - (1-\theta)z}\right)^k,
\end{aligned}$$

and the domain of convergence reduces to $|z| < 1/(1-\theta)$. This is the probability generating function for the $\text{NegBin}(k, 1-\theta)$ distribution. □

**PROOF OF THEOREM 4A:** Here we prove the Maclaurin expansion for the cumulant generating function. Starting with the form in Theorem 3, we have:

$$K(s) = k\ln(\theta) - \sum_{\ell=m-k+1}^{m} \ln\left(\frac{m - me^s + \ell\theta e^s}{\ell}\right)$$

$$= k\ln(\theta) - \sum_{\ell=m-k+1}^{m} \ln\left(\theta e^s + \frac{m}{\ell}\cdot(1-e^s)\right)$$

$$= k\ln(\theta) - \sum_{\ell=m-k+1}^{m} \left[\ln(\theta e^s) + \ln\left(1 - \frac{m}{\ell\theta}\cdot(1-e^{-s})\right)\right]$$

$$= -ks - \sum_{\ell=m-k+1}^{m} \ln\left(1 - \frac{m}{\ell\theta}\cdot(1-e^{-s})\right).$$

Thus, using the well-known Maclaurin expansion $\ln(1-x) = -\sum_{n=1}^{\infty} x^n/n$, the Maclaurin series for the cumulant function is:

$$K(s) = -ks + \sum_{\ell=m-k+1}^{m} \sum_{n=1}^{\infty} \frac{1}{n}\left(\frac{m}{\ell\theta}\cdot(1-e^{-s})\right)^n$$

$$= -ks + \sum_{n=1}^{\infty} \frac{(1-e^{-s})^n}{n} \sum_{\ell=m-k+1}^{m} \left(\frac{m}{\ell\theta}\right)^n$$

$$= -ks + \sum_{n=1}^{\infty} \frac{(1-e^{-s})^n}{n}\left(\frac{m}{\theta}\right)^n (H_m^{(n)} - H_{m-k}^{(n)}).$$

The radius of convergence for this expansion requires $|1 - e^{-s}| < \ell\theta/m$ for all index values $\ell = m - k + 1, \ldots m$. Since $\theta > 0$ and $\ell > 0$ we have $\ell\theta/m > 0$ so the radius of convergence includes an open neighbourhood around $s = 0$. ∎

**LEMMA 1:** Define the function:

$$F(s) = \sum_{n=1}^{\infty} \frac{(1-e^{-s})^n}{n} \cdot f(n).$$

For $r > 0$ the derivatives of this function are:

$$F^{(r)}(s) \equiv \frac{d^r F}{ds^r}(s) = \sum_{i=1}^{r}(-1)^{r-i}\cdot S(r,i)\cdot e^{-is} \sum_{n=i}^{\infty}(n-1)_{i-1}\cdot(1-e^{-s})^{n-i}\cdot f(n).$$

We therefore have:

$$F^{(r)}(0) = \sum_{i=1}^{r}(-1)^{r-i}\cdot S(r,i)\cdot (i-1)!\cdot f(i).$$

**PROOF OF LEMMA 1:** To obtain the result, we will first prove that the derivatives are:

$$F^{(r)}(s) = \sum_{i=1}^{r}(-1)^{r-i} \cdot S(r,i) \cdot e^{-is} \sum_{n=i}^{\infty}(n-1)_{i-1} \cdot (1-e^{-s})^{n-i} \cdot f(n).$$

We prove this result by weak induction. For the base case $r = 1$ we have:

$$F^{(1)}(s) = \sum_{n=1}^{\infty} \frac{d}{ds} \frac{(1-e^{-s})^n}{n} \cdot f(n)$$

$$= e^{-s} \sum_{n=1}^{\infty}(1-e^{-s})^{n-1} \cdot f(n),$$

which is the required equation. For the inductive step we assume that the stipulated form holds for some $r$ and we then have:

$$F^{(r+1)}(s) = \frac{d}{ds} \sum_{i=1}^{r}(-1)^{r-i} \cdot S(r,i) \cdot e^{-i} \sum_{n=i}^{\infty}(n-1)_{i-1} \cdot (1-e^{-s})^{n-i} \cdot f(n)$$

$$= \sum_{i=1}^{r}(-1)^{r-i}(-i \cdot S(r,i))e^{-is} \sum_{n=i}^{\infty}(n-1)_{i-1} \cdot (1-e^{-s})^{n-i} \cdot f(n)$$

$$+ \sum_{i=1}^{r}(-1)^{r-i}S(r,i)e^{-(i+1)s} \sum_{n=i+1}^{\infty}(n-1)_{i} \cdot (1-e^{-s})^{n-i-1} \cdot f(n)$$

$$= \sum_{i=1}^{r}(-1)^{r-i+1}(i \cdot S(r,i))e^{-i} \sum_{n=i}^{\infty}(n-1)_{i-1} \cdot (1-e^{-s})^{n-i} \cdot f(n)$$

$$+ \sum_{i=2}^{r+1}(-1)^{r-i+1}S(r,i-1)e^{-is} \sum_{n=i}^{\infty}(n-1)_{i-1} \cdot (1-e^{-s})^{n-i} \cdot f(n)$$

$$= \sum_{i=1}^{r+1}(-1)^{r-i+1}(i \cdot S(r,i) + S(r,i-1))e^{-is} \sum_{n=i}^{\infty}(n-1)_{i-1} \cdot (1-e^{-s})^{n-i} \cdot f(n)$$

$$= \sum_{i=1}^{r+1}(-1)^{r-i+1}S(r+1,i)e^{-is} \sum_{n=i}^{\infty}(n-1)_{i-1} \cdot (1-e^{-s})^{n-i} \cdot f(n),$$

where the last step follows from the recursive properties of the Stirling numbers of the second kind. This establishes the stipulated form for $r+1$ which establishes the inductive step and thereby proves the form of the derivative. Substitution of $s = 0$ gives the subsequent equation in the lemma. ∎

**PROOF OF THEOREM 4B:** Using the Maclaurin expansion for the cumulant generating function and the functional form in Lemma 1, we have:

$$K(s) = F(s) - ks \qquad f(n) = \left(\frac{m}{\theta}\right)^n (H_m^{(n)} - H_{m-k}^{(n)}).$$

Thus, applying Lemma 1 gives:

$$\kappa_r = F^{(r)}(0) - k \cdot \mathbb{I}(r=1)$$

$$= \sum_{i=1}^{r} (-1)^{r-i} \cdot S(r,i) \cdot (i-1)! \cdot f(i) - k \cdot \mathbb{I}(r=1)$$

$$= \sum_{i=1}^{r} (-1)^{r-i} \cdot S(r,i) \cdot (i-1)! \cdot \left(\frac{m}{\theta}\right)^i (H_m^{(i)} - H_{m-k}^{(i)}) - k \cdot \mathbb{I}(r=1),$$

which was to be shown. ∎

**PROOF OF THEOREM 5:** To facilitate our analysis, we will define the points $x_\ell = \ell/k$ and use the corresponding objects:

$$\Delta x_\ell = \frac{1}{k} \qquad f(x) = \ln\left(\frac{m - kx}{m - m(1-\theta)e^s - kx\theta e^s}\right).$$

We can now write the cumulant generating function in a form that uses a Riemann sum:

$$K(s) = k\left[\ln(\theta) + \sum_{\ell=0}^{k-1} f(x_\ell)\Delta x_\ell\right].$$

Taking the limits $m \to \infty$ and $k \to \infty$ with $k/m \to \lambda$ gives $\Delta x_\ell \to 0$ for each $\ell = 0, \ldots, k-1$ which gives the asymptotic equivalence:

$$K(s) \sim m\lambda\left[\ln(\theta) + \int_0^1 f(x)dx\right].$$

To complete the proof we need to compute the definite integral in this expression. With a bit of calculus, we obtain the indefinite integral:

$$\int f(x)dx = \frac{m}{k}\left[\frac{1-(1-\theta)e^s}{\theta e^s} \cdot \ln|m - m(1-\theta)e^s - kx\theta e^s| - \ln|m - kx|\right]$$

$$+ x \ln\left|\frac{m - kx}{m - m(1-\theta)e^s - kx\theta e^s}\right| + \text{const.}$$

Definite integration between $x = 0$ and $x = 1$ yields the asymptotic form in the theorem. The result for the classical case follows directly by substitution. ∎

**PROOF OF THEOREM 6:** For brevity, we will use the operator lim to refer to the limit under the conditions of the theorem (i.e., with $m \to \infty$ and $k \to \infty$ with the ratio $k/m \to \lambda$). Theorem 6 gives the asymptotic form for the cumulant function, which we will denote as follows:

$$\vec{K}(s) \equiv \lim K(s) = m \begin{bmatrix} \lambda \ln(\theta) - (1-\lambda)\ln|1-\lambda| \\ -\dfrac{1-(1-\theta)e^s}{\theta e^s} \cdot \ln|1-(1-\theta)e^s| \\ +\dfrac{1-(1-\theta)e^s - \lambda\theta e^s}{\theta e^s} \cdot \ln|1-(1-\theta)e^s - \lambda\theta e^s| \end{bmatrix}.$$

Assuming we can interchange the limiting operator and the derivative operator, we obtain:

$$K^{(r)}(s) = \lim \frac{d^r}{ds^r} K(s) = \frac{d^r}{ds^r} \lim K(s) = \frac{d^r}{ds^r} \vec{K}(s) = \vec{K}^{(r)}(s).$$

Substituting $s = 0$ then gives:

$$\lim \kappa_r = K^{(r)}(0) = \vec{K}^{(r)}(0) = \kappa_r^*.$$

This gives the asymptotic equivalence $\kappa_r^* \sim \kappa_r$ so we take $\mu_{m,k,\theta}^* = \kappa_1^*$ and $\sigma_{m,k,\theta}^{*2} = \kappa_2^*$. Now, to establish that we can interchange the limit operator with the derivative operator, we note that for all $m$ and $k$ the cumulant function $K$ is infinitely differentiable at the point $s = 0$, and the sequence $\vec{K}^{(r)}(s)$ is uniformly convergent in a neighbourhood of $s = 0$. Applying the Moore-Osgood limit-interchange theorem then allows us to interchange the limit and the derivative operators. It remains only to derive the relevant cumulants $\kappa_r^*$ and establish the asymptotic distribution. To facilitate the remainder of the proof, we define the function:

$$F(s) \equiv \frac{1 - ae^s - be^s}{be^s} \cdot \ln|1 - ae^s - be^s| - \frac{1 - ae^s}{be^s} \cdot \ln|1 - ae^s|.$$

Taking $a = 1 - \theta$ and $b = \lambda\theta$ we have $\vec{K}(s) = \text{const} + m\lambda F(s)$ which allows us to write the asymptotic cumulants as $\kappa_r^* = \vec{K}^{(r)}(0) = m\lambda F^{(r)}(0)$ for $r \geq 1$. This establishes the general form of the cumulants asserted in the theorem, and it also means that the cumulants of $T_k/m$ are $\kappa_r^{**} = O(m^{1-r})$, which means that they vanish in the limit for $r > 1$. This is sufficient to establish that the random variable $T_k/m$ convergence in distribution to the lognormal, so it remains only to derive the exact form of the mean and variance. With a bit of calculus and algebra, it can be shown that:

$$F^{(1)}(s) = \frac{1}{be^s} \cdot \ln\left|\frac{1-ae^s}{1-ae^s-be^s}\right| - 1,$$

$$F^{(2)}(s) = \frac{1}{(1-ae^s)(1-ae^s-be^s)} - F^{(1)}(s) - 1,$$

$$F^{(3)}(s) = \frac{2ae^s + be^s - 2ae^s(ae^s+be^s)}{(1-ae^s)^2(1-ae^s-be^s)^2} - F^{(2)}(s),$$

$$F^{(4)}(s) = \frac{\begin{bmatrix}(2ae^s+be^s)(1-3ae^s(ae^s+be^s))\\+(be^s)^2+4(ae^s)^2(ae^s+be^s)^2\end{bmatrix}}{(1-ae^s)^3(1-ae^s-be^s)^3} - F^{(3)}(s).$$

Substituting $s=0$ with $a=1-\theta$ and $b=\lambda\theta$ gives:

$$F^{(1)}(0) = -\frac{1}{\lambda\theta} \cdot \ln(1-\lambda) - 1,$$

$$F^{(2)}(0) = \frac{1}{\theta^2(1-\lambda)} + \frac{1}{\lambda\theta} \cdot \ln(1-\lambda),$$

$$F^{(3)}(0) = \frac{2 - 3\theta - \lambda + 3\lambda\theta}{\theta^3(1-\lambda)^2} - \frac{1}{\lambda\theta} \cdot \ln(1-\lambda),$$

$$F^{(4)}(0) = \frac{\begin{bmatrix}6-6\lambda+2\lambda^2\\-12\theta+18\theta\lambda-6\theta\lambda^2\\+7\theta^2-14\theta^2\lambda+7\theta^2\lambda^2\end{bmatrix}}{\theta^4(1-\lambda)^3} + \frac{1}{\lambda\theta} \cdot \ln(1-\lambda).$$

We therefore have the asymptotic moments:

$$\mu^*_{m,k,\theta} \equiv \kappa^*_1 = m\lambda F^{(1)}(0) = \frac{m\lambda}{\theta}\left(-\theta - \frac{\ln(1-\lambda)}{\lambda}\right),$$

$$\sigma^{*2}_{m,k,\theta} \equiv \kappa^*_2 = m\lambda F^{(2)}(0) = \frac{m\lambda}{\theta}\left(\frac{1}{\theta(1-\lambda)} + \frac{\ln(1-\lambda)}{\lambda}\right),$$

$$\kappa^*_3 = m\lambda F^{(3)}(0) = \frac{m\lambda}{\theta}\left(\frac{2-\lambda-3\theta+3\lambda\theta}{\theta^2(1-\lambda)^2} - \frac{\ln(1-\lambda)}{\lambda}\right),$$

$$\kappa^*_4 = m\lambda F^{(4)}(0) = \frac{m\lambda}{\theta}\left(\frac{\begin{bmatrix}6-6\lambda+2\lambda^2\\-12\theta+18\lambda\theta-6\lambda^2\theta\\+7\theta^2-14\lambda\theta^2+7\lambda^2\theta^2\end{bmatrix}}{\theta^3(1-\lambda)^3} + \frac{\ln(1-\lambda)}{\lambda}\right).$$

Substituting $\lambda = k/m$ and simplifying these formulae we obtain the asymptotic forms shown in the theorem. ∎

**PROOF OF THEOREM 7:** To obtain the desired equations in log-space, we will first establish the corresponding recursive equation in regular space. The first equation for the base case follows from the fact that $\text{NegOcc}(t|m, 1, \theta) = \text{Geom}(t|\theta)$, established in Theorem 1. To get the main recursive equation in the theorem we will make use of the "telescoping rule" for the noncentral Stirling numbers of the second kind (O'Neill 2021a, Lemma 2):

$$S(n+1, k, \phi) = \sum_{r=0}^{n-k+1} (k+\phi)^r \cdot S(n-r, k-1, \phi).$$

Applying this rule we obtain the recursive equation:

$$\text{NegOcc}(t|m, k+1, \theta) = \frac{\theta^{k+t+1}}{m^{k+t+1}} \cdot (m)_{k+1} \cdot S\left(k+t, k, m \cdot \frac{1-\theta}{\theta}\right)$$

$$= \frac{\theta^{k+t+1}}{m^{k+t+1}} \cdot (m)_{k+1} \cdot \sum_{r=0}^{t} \left(k + m \cdot \frac{1-\theta}{\theta}\right)^r \cdot S\left(k+t-r-1, k-1, m \cdot \frac{1-\theta}{\theta}\right)$$

$$= \theta \cdot \frac{m-k}{m} \cdot \sum_{r=0}^{t} \left(k + m \cdot \frac{1-\theta}{\theta}\right)^r \cdot \left(\frac{\theta}{m}\right)^r \cdot \text{NegOcc}(t-r|m, k, \theta)$$

$$= \theta \cdot \frac{m-k}{m} \cdot \sum_{r=0}^{t} \left(1 - \theta \cdot \frac{m-k}{m}\right)^r \cdot \text{NegOcc}(t-r|m, k, \theta).$$

Using the notation in the theorem, we can write this more succinctly as:

$$\exp(L(t, k+1)) = \theta \cdot \frac{m-k}{m} \cdot \sum_{r=0}^{t} \exp(r \cdot L_k + L(t-r, k)).$$

Taking logarithms of both sides of this equation gives:

$$L(t, k+1) = \log\left(\theta \cdot \frac{m-k}{m}\right) + \log\left(\sum_{r=0}^{t} \exp(r \cdot L_k + L(t-r, k))\right)$$

$$= \log\left(\theta \cdot \frac{m-k}{m}\right) + \text{logsumexp}\begin{pmatrix} L(t, k), \\ L_k + L(t-1, k), \\ 2 \cdot L_k + L(t-2, k), \\ \vdots \\ t \cdot L_k + L(0, k) \end{pmatrix}.$$

This establishes the main equation in the theorem, which completes the proof. ∎